\documentstyle[12pt,leqno,amssymb,amscd,graphics,hyperref]{amsart}  
\hypersetup{hidelinks=true}

\newtheorem{thm}{Theorem}[section]

\newtheorem{prop}[thm]{Proposition}

\theoremstyle{definition}
\newtheorem{dfn}[thm]{Definition}
\theoremstyle{remark}

\begin{document}

\newcommand{\ct}{\cite}
\newcommand{\pr}{\protect\ref}
\newcommand{\su}{\subseteq}

\newcounter{numb}

\title{Blotto Games with Costly Winnings}

\author{Irit Nowik}
\address{Department of Industrial Engineering and Management, Lev Academic Center, 
P.O.B 16031, Jerusalem 9116001, Israel}
\email{nowik@@jct.ac.il}

\author{Tahl Nowik}
\address{Department of Mathematics, Bar-Ilan University, 
Ramat-Gan 5290002, Israel}
\email{tahl@@math.biu.ac.il}
\urladdr{\url{www.math.biu.ac.il/~tahl}}

\date{March 3, 2016}
\begin{abstract}

We introduce a new variation of the $m$-player asymmetric Colonel Blotto game, where the $n$ battles occur as sequential stages of the game, and the winner of each stage needs to spend resources for
maintaining his win. The limited resources of the players are thus needed both for increasing the probability of winning and for the maintenance costs.
We show that if the initial resources of the players are not too small,
then the game has a unique Nash equilibrium,
and the given equilibrium strategies guarantee the given expected payoff for each player.
\end{abstract}

\maketitle

\section{Introduction}\label{intro}

We present a new $n$-stage game, which is a variation of the Colonel Blotto game.
Each player starts the game with some given resource, and at the beginning of
each stage he must decide how much resource to invest in that stage. A player wins the given stage with probability corresponding to the relative investments of the players, and if all players invest 0 then no player wins that stage.  The winner of the stage receives a payoff
which may differ from stage to stage. Since it is possible that certain stages will not be won by any player,  this is not a fixed sum game.

The players' resources from which the investments are taken can be thought of as money, whereas the payoffs should be thought of as a quantity of different nature, such as political gain. The two quantities cannot be interchanged, that is, the payoff cannot be converted into resources for further investment.

The new feature of our game is the following. The winner of each stage is required to spend additional resources on the maintenance of his winning.
This is a real life situation, where the winnings are some assets, and resources are required for their maintenance, as in wars, territorial contests among organisms, or in the political arena. 
The winner of a given stage must put aside all resources that will be required for future maintenance costs of the won asset.
Thus, a fixed amount will be deducted from the resources of the winner immediately after winning,
which should be thought of as the sum of all future maintenance costs for the given acquired asset. 

At each stage the player thus needs to decide how much to invest in the given stage, where  
winning that stage on one hand leads to the payoff of the given stage, but on the other hand
the maintenance cost for the given winning negatively affects the probabilities for future winnings.
In the present work we show that if the initial resources of the players are not too small then the game has a unique Nash equilibrium, and each player guarantees the payoff of this Nash equilibrium (Theorem~\pr{t1} for two players and Theorem~\pr{tm} for $m$ players.)

As mentioned, our game presents a variation of the well known Colonel Blotto game (\cite{b}).  In Blotto games two players simultaneously distribute forces across several battlefields. 
At each battlefield, the player that allocates the largest force wins. The Blotto game has been developed and generalized in many directions (see e.g., \cite{b}, \cite{f}, \cite{l}, \cite{r}, \cite{h},  \cite{dm}).
Two main developments are the ``asymmetric'' and the ``stochastic'' models. The asymmetric version allows the payoffs of the battlefields to differ from each other, and in the stochastic model the deterministic rule deciding on the winner is replaced by a probabilistic one, by which the chances of winning a battlefield depends on the size of investment.

The present work adds a new feature which changes the nature of the game, in making the winnings costly. 
The players thus do not know before hand how much of their resources will be available for investing in
winning rather than on maintenance, and so the game cannot be formulated with simultaneous investments, as in the usual Blotto games, but rather must be formulated with sequential stages.
At each stage the players need to decide how much to invest in the given stage, based on their remaining available resources and on the future fees and payoffs. 

This work was inspired by previous work of the first author with S. Zamir and I. Segev (\cite{n},\cite{nzs}) on a developmental competition that occurs in the nervous system, which we now describe. A muscle is composed of many muscle-fibers. At birth each muscle-fiber is innervated by several motor-neurons (MNs) that ``compete'' to singly innervate it. It has been found that MNs with higher activation-threshold win in more competitions than MNs with lower activation thresholds. In \cite{n} this competitive process is modeled as a multi-stage game between two groups of players: 
those with lower and those with higher thresholds. At each stage a competition at the most active muscle-fiber is resolved. The strategy of a group is defined as the average activity level of its members and the payoff is defined as the sum of their wins. If a MN wins (i.e., singly innervates) a muscle-fiber, then from that stage on, it must continually devote resources for maintaining this muscle-fiber. Hence the MNs use their resources both for winning competitions and for maintaining previously acquired muscle-fibers. It is proved in \cite{n} that in such circumstances it is advantageous to win in later competitions rather than in earlier ones, since winning at a late stage will encounter less maintenance and thus will negatively affect only the few competitions that were not yet resolved. If $\mu$ is the cost of maintaining a win at each subsequent stage,  then in the terminology of the present work, the fee payed by the MNs for winning the $k$th stage of an $n$  stage game is $(n-k)\mu$.

We start by analyzing the 2-person game in Sections \pr{tg} and \pr{ns}, and then generalize to the $m$-person setting in Section~\pr{g}.

\section{The 2-person game}\label{tg}

The initial data for the 2-person version of our game is the following.
\begin{enumerate}
\item The number $n$ of stages of the game.
\item Fixed payoffs $w_k>0$, \ $1\leq k\leq n$, to be received by the winner of the $k$th stage.
\item The initial resources $A, B \geq0$ of players I,II respectively. 
\item Fixed fees $c_k\geq 0$, \ $1\leq k \leq n-1$, to be deducted from the resources of the winner after the $k$th stage.

\end{enumerate}

The rules of the game are as follows. At the $k$th stage of the game, the two 
players, which we  name PI,PII, each has some remaining resource $A_k,B_k$, where $A_1=A,B_1=B$.
PI,PII each needs to decide his investment $x_k,y_k$ for that stage, respectively, with $0\leq x_ k\leq A_k -c_k$, \ $0\leq y_ k\leq B_k -c_k$, and where if $A_k < c_k$ then PI may only invest 0, and similarly for PII. 
These rules ensure that the winner of the given stage will have the resources for paying the given fee $c_k$.
The probability for PI,PII of winning this stage is respectively $\frac{x_k}{x_k+y_k}$ and $\frac{y_k}{x_k+y_k}$, where if  $x_k=y_k=0$ then no player wins. 
The resource of the winner of the $k$th stage is then
reduced by an additional $c_k$, that is, if PI wins the $k$th stage then $A_{k+1}=A_k-x_k-c_k$ and $B_{k+1}=B_k-y_k$, and if PII wins then $A_{k+1}=A_k-x_k$ and $B_{k+1}=B_k-y_k-c_k$. 
The role of $c_k$ is in determining $A_{k+1},B_{k+1}$, thus there is no $c_n$.  
It will however be convenient in the sequel to formally define $c_n=0$.
The payoff received by the winner of the $k$th stage is $w_k$. 
Since it is possible that no player wins certain stages, this game is not a fixed sum game.

As already mentioned, the resource quantities $A_k,B_k,x_k,y_k,c_k$ used for the investments and fees are of different 
nature than that of the payoffs $w_k$. The two quantities cannot be interchanged and should be thought of as having different ``units''. 
Note that all expressions below are unit consistent, that is, if say we divide resources by payoff, then such expression has units of $\frac{\text{resources}}{\text{payoff}}$, and may only be added or equated to expressions of the same units.

If $A$ and $B$ are too small in comparison to $c_1,\dots,c_{n-1}$ then
the players' strategies are strongly influenced by the possibility of running out of
resources before the end of the game. In the present work we analyze the game when $A,B$ are not too small. Namely, we
introduce a quantity $M$ depending on $c_1,\dots,c_{n-1}$ and $w_1,\dots,w_n$, and prove that if
$A,B>M$ then there is a unique Nash equilibrium for the game, and each player guarantees the value of this Nash equilibrium.

For $k=1,\dots,n$ let $W_k=\sum_{i=k}^n w_i$ and $W=W_1$.
We now show that if $A>M$, then if PI always chooses to invest $x_k \leq\frac{w_k}{W_k} A_k$ (as holds for our strategy $\sigma_{n,A,B}$ presented in Definition~\pr{st} below), then whatever the random outcomes of the game are, 
his resources will not run out before the end of the game. We in fact give a specific lower bound on $A_k$ for every $k$, which will be used repeatedly in the sequel.

\begin{prop}\label{m}
Let
$$M =W\cdot\max_{1\leq k \leq n} \Bigg(\frac{c_k}{w_k} + \sum_{i=1}^{k-1} \frac{ c_i}{W_{i+1}} \Bigg).$$
If $A>M$, and if PI plays $x_k \leq \frac{w_k A_k}{W_k}$ for all $k$, 
then
$A_k > \frac{W_k c_k}{w_k}$
for all $1 \leq k \leq n$. In particular $A_k>0$ 
for all $1 \leq k \leq n$. 
And similarly for PII. 
\end{prop} 

\begin{pf} For every $1 \leq k \leq n$ 
we have $\frac{A}{W} > \frac{M}{W} \geq \frac{c_k}{w_k}+\sum_{i=1}^{k-1} \frac{c_i}{W_{i+1}}$, 
so $$\frac{A}{W} - \sum_{i=1}^{k-1} \frac{c_i}{W_{i+1}} > \frac{c_k}{w_k}.$$ Thus it is enough to show that $\frac{A_k}{W_k} \geq \frac{A}{W} - \sum_{i=1}^{k-1} \frac{c_i}{W_{i+1}}$
for all $1 \leq k \leq n$.
We show this by induction on $k$. For $k=1$ the sum is empty and we get equality.
Assuming $$\frac{A_k}{W_k} \geq \frac{ A}{W} - \sum_{i=1}^{k-1} \frac{c_i}{W_{i+1}}$$ we get
\begin{multline*}
\frac{A_{k+1}}{W_{k+1}}
\geq \frac{1}{W_{k+1}}\Big(A_k - \frac{w_kA_k}{W_k} -c_k\Big) 
=  \frac{1}{W_{k+1}}\Big(\frac{W_{k+1}A_k}{W_k} - c_k\Big)    \\  
=\frac{A_k}{W_k}-\frac{c_k}{W_{k+1}}  
\geq \frac{ A}{W} - \sum_{i=1}^{k-1} \frac{c_i}{W_{i+1}}-\frac{c_k}{W_{k+1}}
=\frac{ A}{W} - \sum_{i=1}^k \frac{c_i}{W_{i+1}}.\end{multline*}
\end{pf}

We consider two simple examples of $c_1,\dots,c_{n-1}$, $w_1,\dots,w_n$, for which 
$M$ may be readily identified.
\begin{enumerate}
\item Let  $c_k=n-k$, \ $w_k=1$ for all $k$. 
These fees and payoffs are as in the biological game described in the introduction.
For every $1\leq k \leq n$
we have
$n\big(c_k+\sum_{i=1}^{k-1} \frac{c_i}{n-i} \big)   = n\big(  (n-k) +  (k-1) \big) = n(n-1)$, so $M=n(n-1)$. 
\item Let $c_k=1$ for all $1 \leq k \leq n-1$, \ $w_k=1$ for all $k$.
Using the inequality $\sum_{k=1}^{n} \frac{1}{k}< 1+ \ln n$, we get for every $1\leq k \leq n$ that
$n\big(c_k+\sum_{i=1}^{k-1} \frac{c_i}{n-i}  \big)  < n(2+ \ln n)$, so $M<n(2+ \ln n)$.
\end{enumerate}

We note that an obvious necessary condition for $A$ to satisfy Proposition~\pr{m} is $A\geq \sum_{k=1}^{n-1} c_k$, since in case PI wins all stages he will need to pay all fees $c_k$.
We see that $M$ in the two examples above is not much larger. Namely, in~(1),
$\sum_{k=1}^{n-1} c_k=\frac{n(n-1)}{2}$ and $M=n(n-1)$, and in~(2), $\sum_{k=1}^{n-1} c_k=n-1$ and $M<n(2+ \ln n)$.

\section{Nash equilibrium}\label{ns}

We define the following two strategies $\sigma_{n,A,B}$ and $\tau_{n,A,B}$ for  PI,PII respectively. We prove that for $A,B>M$ as given in Proposition~\pr{m},
 this pair of strategies is a unique Nash equilibrium, and these strategies guarantee the
given payoffs. 

\begin{dfn}\label{st}
At the $k$th stage of the game, let
$$a_k=\frac{w_kA_k}{W_k} - \frac{A_kc_k}{A_k+B_k}   \ \ \ \ \ \hbox{and} 
\ \ \ \ \ b_k=\frac{w_kB_k}{W_k} - \frac{B_kc_k}{A_k+B_k}.$$ 
where as mentioned, we formally define $c_n=0$.
The strategy $\sigma_{n,A,B}$ for PI is the following: At the $k$th stage PI invests $a_k$ if it is allowed by the rules of the game. Otherwise he invest 0. The strategy $\tau_{n,A,B}$ for PII is similarly defined with $b_k$.
\end{dfn}

Recall that $a_k \neq 0$ is allowed by the rules of the game if $0 \leq a_k \leq A_k-c_k$, 
whereas $a_k=0$ is always allowed, even when $A_k-c_k<0$.
We interpret the quantities $a_k,b_k$ as follows. PI first divides his remaining resource $A_k$ to the remaining stages
in proportion to the payoff for each remaining stage, which gives $\frac{w_k}{W_k}A_k$. 
From this he subtracts $\frac{A_k}{A_k+B_k}c_k$  which is the expected fee he will pay for this stage, since
$\frac{a_k}{a_k+b_k} = \frac{A_k}{A_k+B_k}$.
Note that $W_n=w_n$ and formally $c_n=0$, so $a_n=A_n$, \ $b_n=B_n$, i.e. at the last stage the two players invest all their remaining resources. 

Depending on $A$ and $B$ and on the random outcomes of the game, it may be that PI indeed reaches a stage where $a_k$ is not allowed. In this regard we make the
following definition.

\begin{dfn}\label{d}
The triple $(n,A,B)$ is \emph{PI-effective} if when PI and PII use $\sigma_{n,A,B}$ and $\tau_{n,A,B}$, then it is impossible that they reach a stage where $a_k$ is not allowed for PI.
 Similarly PII-effectiveness is defined for PII with $b_k$.
\end{dfn}

\begin{prop}\label{b2}
Let $M$ be as in Proposition~\pr{m}.
If $A >M$ and $B$ is arbitrary, then $(n,A,B)$ is PI-effective. Furthermore, $a_k>0$ for all $k$.  And similarly for PII when $B>M$.
\end{prop}

\begin{pf}
We need to show that necessarily $0 < a_k \leq A_k-c_k$ for all $1\leq k \leq n$.
We have $a_k=\frac{w_kA_k}{W_k} - \frac{A_kc_k}{A_k+B_k} \leq  \frac{w_kA_k}{W_k}$, so
by Proposition~\pr{m},   $\frac{w_k}{W_k}>\frac{c_k}{A_k}  \geq\frac{c_k}{A_k+B_k}$ and $A_k>0$, so $\frac{w_kA_k}{W_k}>\frac{A_kc_k}{A_k+B_k}$ giving $a_k>0$.

For the inequality $a_k \leq A_k-c_k$ we first consider $k\leq n-1$. We have
from the proof of Proposition~\pr{m} that 
$\frac{A_k}{W_k}-\frac{c_k}{W_{k+1}} \geq \frac{ A}{W} - \sum_{i=1}^k \frac{c_i}{W_{i+1}}
> \frac{c_{k+1}}{w_{k+1}}\geq 0$, so  $\frac{A_k}{W_k}>\frac{c_k}{W_{k+1}}$, and so 
$$(1-\frac{A_k}{A_k+B_k})c_k \leq c_k  <\frac{W_{k+1}A_k}{W_k} = (1-\frac{w_k}{W_k})A_k.$$
This gives $c_k- \frac{A_kc_k}{A_k+B_k} < A_k - \frac{w_kA_k}{W_k}$,
so $a_k=\frac{w_kA_k}{W_k} - \frac{A_kc_k}{A_k+B_k} < A_k-c_k$.
For $k=n$ we note that 
$c_n=0$ by definition, and $W_n=w_n$, so $a_n = A_n =A_n-c_n$. 
\end{pf}

In general, an inductive characterization of PI-effectiveness will also involve induction regrading PII. But if we assume that $B>M$, and so by Proposition~\pr{b2} all $b_k$ are known to be allowed and positive, 
then the notion of PI-effectiveness becomes simpler, and may be characterized inductively as follows.
When saying that a triple $(n-1,A',B')$ is PI-effective, we refer to the $n-1$ stage game 
with fees $c_2,\dots,c_{n-1}$ and payoffs $w_2,\dots,w_n$.
Starting with $n=1$, \ $(1,A,B)$ is always PI-effective. For $n\geq 2$, 
if $a_1$ is not allowed then $(n,A,B)$ is not PI-effective. 
If $a_1=0$ then it is allowed, and PI surely loses the first stage, and so $(n,A,B)$ is PI-effective iff $(n-1,A,B-b_1-c_1)$ is PI-effective. 
Finally if $a_1>0$ and it is allowed then $(n,A,B)$ is PI-effective iff both $(n-1,A-a_1-c_1,B-b_1)$ and 
$(n-1,A-a_1,B-b_1-c_1)$ are PI-effective.

The crucial step in proving  Theorem~\pr{t1} below, on the unique Nash equilibrium and the guaranteed payoffs, is the following Theorem~\pr{bm}. We point out that in Theorem~\pr{t1} we will assume that $A>M$, in which case $(n,A,B)$ is PI-effective, by Proposition~\pr{b2}. But here in Theorem~\pr{bm} we must consider arbitrary $A\geq 0$ in order for an induction argument to carry through.

\begin{thm}\label{bm}
Given $c_1,\dots,c_{n-1}$ and $w_1,\dots,w_n$ let $M$ be as in Proposition~\pr{m}, and assume that $B>M$ and PII plays the strategy $\tau_{n,A,B}$. For $A\geq 0$, if $(n,A,B)$ is PI-effective, and PI plays according to $\sigma_{n,A,B}$, then his expected payoff is $\frac{AW}{A+B}$. On the other hand, if $(n,A,B)$ is not PI-effective, or if PI uses a different strategy, then his expected payoff is strictly less than $\frac{AW}{A+B}$.
\end{thm}

\begin{pf} 
By induction on $n$. We note that throughout the present proof we do not use the condition $B>M$ directly, but rather only through the statements of Propositions \pr{b2} and \pr{m} saying that $(n,A,B)$ is PII-effective, 
$b_k>0$ and $c_k <  \frac{w_k B_k}{W_k}$ for all $1 \leq k \leq n$, which indeed continue to hold along the induction process. 

If $A=0$ then $a_k=0$ for all $k$, which is the only possible investment, and its payoff is $0=\frac{AW}{A+B}$, so the statement holds. 
We thus assume from now on that $A>0$. 
For $n=1$ we have $b_1=B$. The allowed investment for PI is $0\leq s\leq A$ with expected payoff $\frac{s}{s+B}w_1=\frac{s}{s+B}W$ which indeed attains a strict maximum $\frac{A}{A+B}W$ at $s=A=a_1$.

For $n\geq 2$, let $s$ be the investment of PI in the first stage.  Assume first that $s=0$.
In this case
PII surely wins the first stage and so following this stage we have $A_2=A$ and $B_2=B-b_1-c_1$.
The moves for PII dictated by $\tau_{n,A,B}$ for the remaining $n-1$ stages of the game are  $\tau_{n-1,A,B-b_1-c_1}$, and so
by the induction hypothesis the expected total payoff of PI is at most $\frac{AW_2}{A+B-b_1-c_1}$.
Since Proposition~\pr{m} holds for PII, we have $c_1 < \frac{w_1B}{W}\leq \frac{w_1(A+B)}{W}$, that is, $\frac{w_1}{W} - \frac{c_1}{A+B}>0$,
and since $A>0$ we get $a_1 = A(\frac{w_1}{W} - \frac{c_1}{A+B})> 0$. 
This means that $s=0\neq a_1$, so we must verify the strict inequality 
$\frac{AW_2}{A+B-b_1-c_1}<\frac{AW}{A+B}$. 
This is readily verified, using $A>0$, \ $c_1 <  \frac{w_1B}{W}$, \ $W_2=W-w_1$, and 
$b_1+c_1=\frac{w_1B}{W} - \frac{Bc_1}{A+B}+c_1=\frac{w_1B}{W} + \frac{Ac_1}{A+B}$.

We now assume $s>0$. This is allowed only if 
$A > c_1$ and $0 < s \leq A-c_1$. 
The moves for PII dictated by $\tau_{n,A,B}$ for the remaining $n-1$ stages of the game are
$\tau_{n-1,A_2,B_2}$.
By the induction hypothesis,  if PI wins the first stage, which happens with probability $\frac{s}{s+b_1}>0$, then
his expected payoff in the remaining $n-1$ stages of the game is at most 
$\frac{(A-s-c_1)W_2}{A+B-s-b_1-c_1}$. Similarly, if he loses the first stage, which happens 
with probability $\frac{b_1}{s+b_1}>0$, then his expected payoff in the remaining $n-1$ stages
is at most $\frac{(A-s)W_2}{A+B-s-b_1-c_1}$.
Thus, the expected payoff of PI for the whole $n$ stage game is at most $F(s)$, where 
$$F(s)=\frac{s}{s+b_1}\bigg(w_1+\frac{(A-s-c_1)W_2}{A+B-s-b_1-c_1}\bigg) +  \frac{b_1}{s+b_1}\cdot\frac{(A-s)W_2}{A+B-s-b_1-c_1}$$
with $b_1= \frac{w_1B}{W} - \frac{Bc_1}{A+B}$.

By the induction hypothesis we know furthermore, that in case PI wins the first stage, he will attain the maximal expected payoff  $\frac{(A-s-c_1)W_2}{A+B-s-b_1-c_1}$ in the remaining stages of the game only if 
$(n-1,A-s-c_1, B-b_1)$ is PI-effective, and he uses $\sigma_{n-1,A-s-c_1, B-b_1}$. Similarly, if he loses the first stage, he will attain the maximal expected payoff
$\frac{(A-s)W_2}{A+B-s-b_1-c_1}$  
only if $(n-1,A-s,B-b_1-c_1)$ is PI-effective and he uses $\sigma_{n-1,A-s,B-b_1-c_1}$. If not,
then since both alternatives occur with positive probability, his expected total payoff for the whole $n$ stage game will be strictly less than $F(s)$.

To analyze $F(s)$,
we make a change of variable $s=a_1+x$, that is, we define 
$\widehat{F}(x) =F(a_1+x)=F(\frac{w_1A}{W} - \frac{Ac_1}{A+B}+x)$.  
After some manipulations we get:
$$\widehat{F}(x)=\frac{AW}{A+B}- \frac{BW^3x^2}{(A+B)\Big(W_2(A+B)-Wx\Big)\Big(Wx-Wc_1+w_1(A+B)\Big)}.$$

Under this substitution, $s=a_1$ corresponds to $x=0$, and the allowed domain $0< s\leq A-c_1$ corresponds to
$$ \frac{Ac_1}{A+B} -\frac{w_1A}{W}  < x \leq \frac{W_2A}{W}-\frac{Bc_1}{A+B}.$$ 
Using $c_1<\frac{w_1B}{W}$, one may verify that in the above expression for $\widehat{F}$ 
the two linear factors appearing in the denominator of the second term are both strictly positive in this domain. It follows that $\widehat{F}$ in the given domain is at most
$\frac{AW}{A+B}$, and this maximal value is attained only for $x=0$ (if it is in the domain), which corresponds to $s=a_1$ for the original $F$. Finally, as mentioned, 
unless $(n-1,A_2,B_2)$ is PI-effective and PI plays $\sigma_{n-1,A_2,B_2}$, his expected payoff will be strictly less than $F(s)$, and so strictly less than $\frac{AW}{A+B}$.
\end{pf}
 
We may now prove our main result.

\begin{thm}\label{t1}
Given $c_1,\dots,c_{n-1}$ and $w_1,\dots,w_n$, let $M$ be as in Proposition~\pr{m}, and assume $A,B>M$. Then the pair of strategies $\sigma_{n,A,B},\tau_{n,A,B}$ is a unique Nash equilibrium for the game, with expected total  payoffs $\frac{AW}{A+B}$, \ $\frac{BW}{A+B}$.
Furthermore, $\sigma_{n,A,B}$ and $\sigma_{n,A,B}$ guarantee the expected payoffs $\frac{AW}{A+B}$ and $\frac{BW}{A+B}$.
\end{thm}

\begin{pf}
Denote $\sigma_0=\sigma_{n,A,B}$ and $\tau_0=\tau_{n,A,B}$, and for any pair of strategies 
$\sigma,\tau$ let $S_1(\sigma,\tau)$, $S_2(\sigma,\tau)$ be the expected payoffs of PI, PII respectively.
We first prove the second statement of the theorem. 
Recall that if both players invest 0 in a given stage then there is no winner to that stage. However, if $B>M$ and PII plays $\tau_0$, then by Proposition~\pr{b2} we have $b_k>0$ for all $k$,
and so indeed there is a winner to each stage of the game, and thus the total combined payoff of PI and PII is necessarily $W$. It thus follows from Theorem~\pr{bm} that for any strategy $\sigma$ of PI we have 
$S_2(\sigma,\tau_0) = W-S_1(\sigma,\tau_0) \geq \frac{BW}{A+B}$. 
Similarly, if $A>M$ then $S_1(\sigma_0,\tau) \geq \frac{AW}{A+B}$ for all $\tau$, establishing the second statement of the theorem.

As to the first statement, Theorem~\pr{bm} applied to both PI and PII implies that the pair 
$\sigma_0,\tau_0$ is a Nash equilibrium with the given expected payoffs. 
To show it is unique we argue as follows. Let $\sigma,\tau$ be any other Nash equilibrium and assume that $\sigma \neq \sigma_0$.
By Theorem~\pr{bm} we have $S_1(\sigma,\tau_0)<S_1(\sigma_0,\tau_0)$ and since 
playing $\tau_0$ guarantees a combined total payoff of $W$, 
we have $S_2(\sigma,\tau_0)=W-S_1(\sigma,\tau_0)>W-S_1(\sigma_0,\tau_0)=S_2(\sigma_0,\tau_0)$.
Since the pair $\sigma,\tau$ is a Nash equilibrium we also have $S_2(\sigma,\tau)\geq S_2(\sigma,\tau_0)$,
and together we get $S_2(\sigma,\tau)>S_2(\sigma_0,\tau_0)$. Since $S_1(\sigma,\tau)+S_2(\sigma,\tau)\leq W$ and  $S_1(\sigma_0,\tau_0) +S_2(\sigma_0,\tau_0)=W$, we must have 
$S_1(\sigma,\tau)<S_1(\sigma_0,\tau_0)$. Again since $\sigma,\tau$ is a Nash equilibrium
we have $S_1(\sigma_0,\tau) \leq  S_1(\sigma,\tau)$ so together 
$S_1(\sigma_0,\tau) < S_1(\sigma_0,\tau_0)=\frac{AW}{A+B}$, contradicting  the conclusion of the previous paragraph.
\end{pf}

We remark that our Nash equilibrium is subgame perfect. This may be seen from the inductive proof, and
also follows from the uniqueness of the Nash equilibrium, since the game is finite. 

\section{The $m$-person game}\label{g}

The generalization of our game and results to an $m$-person setting is straightforward.
Player $P_i$, \ $1\leq i\leq m$, starts with resource $A^i\geq0$, and we are again given fixed fees  
$c_1,\dots,c_{n-1} \geq 0$, and payoffs $w_1,\dots,w_n>0$. At stage $1\leq k \leq n$ player $P_i$ has remaining resource $A^i_k$,
with $A^i_1=A^i$. At stage $k$
each player decides to  invest $0\leq x^i_k\leq A^i_k - c_k$, or 0 if $A^i_k < c_k$. 
The probability for $P_i$ to win is $\frac{x^i_k}{\sum_{j=1}^m x^j_k}$, and if $x^i_k=0$ for all $m$ then no player wins. The winner of the $k$th stage receives payoff $w_k$ and pays the maintenance fee $c_k$.

We define $M$ as before, $M =W\cdot\max_{1\leq k \leq n} \Bigg(\frac{c_k}{w_k} + \sum_{i=1}^{k-1} \frac{ c_i}{W_{i+1}} \Bigg)$.

The strategy $\sigma^i_{n,A^1,\dots,A^m}$ of $P_i$ is the straightforward generalization of the strategies 
$\sigma_{n,A,B},\tau_{n,A,B}$
of the 2-person game, namely, at the $k$th stage $P_i$ invests
$$a^i_k=\frac{w_kA^i_k}{W_k} - \frac{A^i_k c_k}{\sum_{j=1}^m A^j_k}$$
if it is allowed, and 0 otherwise.

The generalization of Theorem~\pr{t1} is the following.

\begin{thm}\label{tm}
If  $A^1,\dots,A^m > M$ then
the $m$ strategies $\sigma^1_{n,A^1,\dots,A^m},\dots,\sigma^m_{n,A^1,\dots,A^m}$ are a unique Nash equilibrium for the $m$-person game, with expected total payoffs
$\frac{A^1 W}{\sum_{i=1}^m A^i},\dots,\frac{A^m W}{\sum_{i=1}^m A^i}$.
Furthermore, $\sigma^i_{n,A^1,\dots,A^m}$ guarantees the expected payoff $\frac{A^i W}{\sum_{j=1}^m A^j}$ for $P_i$.
\end{thm}

\begin{pf}
We prove for $P_1$. Let $B_k = \sum_{i=2}^m A^i_k$ and $b_k = \sum_{i=2}^m a^i_k$.
Since $a^1_k=\frac{w_k A^1_k}{W_k} - \frac{A^1_k c_k}{A^1_k+B_k}$
and  $b_k=\frac{w_k B_k}{W_k} - \frac{B_kc_k}{A^1_k+B_k}$, 
our player $P_1$ can imagine that he is playing a 2-person game against one joint player whose resource
is $B_k$, and whose strategy dictates investing $b_k$. 
Thus the claim follows from Theorem~\pr{t1}.
\end{pf}

\end{document}